\documentclass[a4paper,11pt]{amsart} 
\usepackage{amscd}             
\usepackage{amssymb}
\usepackage{amsthm}
\usepackage{epsf} 
\usepackage[T1]{fontenc}            
\newtheorem{thm}{Theorem}[section]   
\newtheorem{cor}[thm]{Corollary}
\newtheorem{lemma}[thm]{Lemma}
\newtheorem{prop}[thm]{Proposition}

\newtheorem{defn}[thm]{Definition}
\newtheorem{rem}[thm]{Remark}
\renewcommand{\proofname}{Proof}

\def\ker{\operatorname{ker}}

\def\min{\operatorname{min}}
\def\im{\operatorname{Im}}
\def\rank{\operatorname{rank}}

\def\length{\operatorname{length}}
\def\c1{\operatorname{c_1}}
\def\c2{\operatorname{c_2}}
\def\Cliff{\operatorname{Cliff}}
\def\gon{\operatorname{gon}}

\def\Grass{\operatorname{Grass}}
\def\red{\operatorname{red}}
\def\CC{{\mathbf C}}
\def\ZZ{{\mathbf Z}}

\def\QQ{{\mathbf Q}}
\def\PP{{\mathbf P}}

\def\F{{\mathcal F}}

\def\O{{\mathcal O}}
\def\I{{\mathcal J}}
\def\Z{{\mathcal Z}}

\def\L{{\mathcal L}}
\def\v{^{\vee}}                  
\def\iso{\simeq}                 
\def\eqv{\equiv}
\def\sub{\subseteq}

\def\+{\oplus}                   
\def\*{\otimes}                  
\def\hpil{\longrightarrow}       
\def\khpil{\rightarrow}

\def\Hom{\operatorname{Hom}}

\def\Pic{\operatorname{Pic}}

\def\det{\operatorname{det}}

\begin{document}

  \title{On $k$th-order embeddings of $K3$ surfaces and Enriques surfaces}
  \author{Andreas Leopold Knutsen}  
  
  \address{Department of Mathematics\\ 
    University of Bergen\\ Johs. Brunsgt 12\\ N-5008 Bergen\\ Norway}

  \begin{abstract}
   We give necessary and sufficient conditions for a big and nef line 
   bundle $L$ of any degree on
   a $K3$ surface or on an Enriques surface $S$ to be $k$-very ample and 
   $k$-spanned. 
   Furthermore, we give necessary and sufficient conditions for a spanned 
   and big line bundle on
   a $K3$ surface $S$ to be birationally $k$-very ample and birationally 
   $k$-spanned (our definition), and relate these concepts
   to the Clifford index and gonality of smooth curves in $|L|$ and the 
   existence of a particular type of rank $2$ bundles on $S$.
\end{abstract}

\maketitle

\section{Introduction}
\label{intro}

Let $L$ be a line bundle on a smooth connected surface $S$ over the complex numbers. Recall \cite{BS2} that $L$ is called \textit{$k$-very ample}, for an integer $k \geq 0$, if for any $0$-dimensional subscheme $( \Z, \O _{\Z} )$ of length $h^0 (\O _{\Z} )= k+1$, the restriction map $H^0 (L) \khpil H^0 (L \* \O _{\Z} )$ is surjective.  This definition is a natural generalization of the notions of spannedness and very-ampleness of line bundles. In fact, by definition, $L$ is $0$-very ample if and only if $L$ is generated by its global sections, and $L$ is $1$-very ample if and only if $L$ is very ample.

There are various geometrical interpretations of the notion of $k$-very ampleness. Denoting by $S^{[r]}$ the Hilbert scheme of $0$-dimensional subschemes of $S$ of length $r$, and by $\Grass (r, H^0 (L))$ the Grassmannian of all $r$-dimensional quotients of $H^0 (L)$, then the rational map
\[ \phi _k : S^{[k+1]} \hpil \Grass (k+1, H^0 (L)), \]
sending $( \Z, \O _{\Z} ) \in S^{[k+1]}$ into the quotient $H^0 (L) \khpil H^0 (L \* \O _{\Z} )$ is a morphism if $L$ is $k$-very ample and an embedding if and only if $L$ is $(k+1)$-very ample \cite{cg}.

Also, if $S$ is embedded in $\PP ^{h^0(L)-1}$ via a very ample line bundle $L$, then $L$ is $k$-very ample if and only if $S$ has no $(k+1)$-secant $(k-1)$-planes.

There is also a slightly weaker condition than $k$-very ampleness as follows \cite{BFS}:
$L$ is called \textit{$k$-spanned}, for an integer $k \geq 0$, if for any \textit{curvilinear} $0$-dimensional subscheme $( \Z, \O _{\Z} )$ of length $h^0 (\O _{\Z} )= k+1$, the restriction map $H^0 (L) \khpil H^0 (L \* \O _{\Z} )$ is surjective. Recall that a $0$-dimensional scheme $( \Z, \O _{\Z} )$ is called curvilinear if $\dim T_{x} \Z \leq 1$ for every $x \in \Z _{\red}$. On a smooth irreducible curve, the two notions of $k$-very ampleness and $k$-spannedness coincide for all $k \geq 0$. On a smooth connected surface they coincide for $k \leq 2$ \cite[Lemma 3.1]{BFS}. In the sequel we will show that for $K3$ and Enriques surfaces these two notions are equivalent for all $k \geq 0$.

In recent years a lot of work has been done in the study of $k$-very
ample and $k$-spanned line bundles on surfaces (see e.g. \cite{ba}, \cite{BFS}, \cite{BS1}, \cite{BS2}, \cite{BS3}, \cite{BS4}, \cite{dr1}, \cite{dr2}, \cite{ter}, \cite{BSz}), and in particular, in the
classification of pairs $(S,L)$, where $S$ is a surface and $L$ is a $k$-very
ample or $k$-spanned line bundles on $S$.

In \cite{BFS} (resp. \cite{BS2}) Beltrametti, Francia and Sommese (resp.
the first and third author) showed that if $L$ is nef and $L^2 \geq
4k+5$, and $K_S +L$ 
is not $k$-very ample (resp. $k$-spanned), then there exists an
effective divisor $D$ such that 
\[ L.D - k-1 \leq D^2 < L.D/2 < k+1. \]

Recently, Terakawa \cite{ter} showed that for line bundles of degree $>4k+4$ 
on surfaces of Kodaira
dimension zero, these conditions
are also sufficient. Since $K_S$ is numerically equivalent to zero, 
in particular necessary and
sufficient conditions for a nef line bundle $L$ to be $k$-very ample and
$k$-spanned were granted. In fact these conditions are equivalent for
the $k$-very ample and the $k$-spanned case.

If $L$ is a big $k$-spanned line bundle on a smooth surface $S$ of Kodaira 
dimension zero and $L^2 \leq 4k+4$, then either $S$ is a $K3$ surface and 
$L^2= 4k$, $4k+2$ or $4k+4$, or $S$ is an Enriques surface and $L^2= 4k+4$ 
(see Proposition \ref{under} below).

In this paper we complete the description of $k$-very ample and
$k$-spanned line bundles on surfaces of Kodaira
dimension zero. Our approach holds for all $L^2 \geq 4k$ on $K3$ surfaces 
(resp. all 
$L^2 \geq 4k+4$ on Enriques surfaces), so no condition that
$L^2 \leq 4k+4$ will be imposed, since they would not make the proofs easier. 
Thus, we give a unified presentation of the cases $L^2 \geq 4k+5$, already 
treated by the mentioned authors, and the cases with low values of $L^2$, 
where the results are new. We will need a different
approach than in \cite{BFS} and \cite{BS2}, where Bogomolov stability is
used, making the assumption $L^2 \geq 4k+5$ necessary. 

The author would like to mention that at the same time that an earlier version of this paper was written, and independently, T.~Szemberg \cite{sz} treated the Enriques case. In particular, he showed that the case $L^2= 4k+4$ and $L$ $k$-very ample only occurs if $k=0$.  

For spanned and big line bundles of any degree on a $K3$ surface our approach 
also makes it possible to give a characterization of birational $k$-very 
ampleness and birational $k$-spannedness. A 
big and globally generated line bundle $L$ will be called 
\textit{birationally $k$-very ample} (resp. \textit{birationally 
$k$-spanned)}, if there exists a non-empty Zariski-open 
subset of $S$ where $L$ is $k$-very ample (resp. $k$-spanned).  
These concepts (which on a $K3$ surface turn   
out to be equivalent as well) are interesting not only in their own rights, 
but also because
they are connected to the Clifford index of 
smooth curves in $|L|$, the minimal gonality of smooth curves in $|L|$ 
\cite{cp} and the existence of a certain type of rank $2$ vector bundle on $S$.
   
The following three theorems are the main results in this paper.

\begin{thm} \label{kvathm}
  Let $L$ be a big and nef line bundle on a $K3$ surface and $k \geq 0$ 
  an integer. The following
  conditions are equivalent:
\begin{itemize}
  \item [a)] $L$ is $k$-very ample,
  \item [b)] $L$ is $k$-spanned,
  \item [c)] $L^2 \geq 4k$ and there exists no effective divisor $D$
    satisfying the conditions $(*)$ below:
     \begin{eqnarray*}
& 2D^2 \stackrel{i)}{\leq} L.D \leq D^2 + k+1 \stackrel{ii)}{\leq} 2k+2   \\
(*) & \mbox{ with equality in i) if and only if } L \sim 2D \mbox{ and } L^2 \leq 4k+4,   \\
& \mbox{and
equality in ii) if and only if } L \sim 2D \mbox{ and } L^2= 4k+4. 
     \end{eqnarray*}
\end{itemize}

Furthermore, if $L$ is not $k$-very ample (equivalently $k$-spanned),
then among all divisors satisfying $(*)$, we can always find a smooth
curve, and all smooth curves satisfying $(*)$
will contain a $0$-dimensional scheme of degree $k+1$ where the 
$k$-spannedness fails (more precisely $L_{|D}$ is not $k$-spanned).
\end{thm}

The corresponding result for Enriques surfaces is the following, which in the 
$k$-very ample case is similar to a result obtained independently by Szemberg 
in \cite{sz}. Our approach is slightly different.

\begin{thm} \label{kvathmenr}
  Let $L$ be a big and nef line bundle on an Enriques surface and $k \geq 0$ 
  an integer. The following conditions are equivalent:
\begin{itemize}
  \item [a)] $L$ is $k$-very ample,
  \item [b)] $L$ is $k$-spanned,
  \item [c)] there exists no divisor $D>0$ satisfying $D^2=-2$, 
             $D.L \leq k-1$, or $D^2=0$, $D.L \leq k+1$.
\end{itemize}
\end{thm}

So, unlike in the $K3$ case, the $k$-very ampleness of $L$ is governed by divisors with self-intersection $0$ and $-2$. However, we will see that divisors satisfying similar conditions as the conditions $(*)$ play an important role also for Enriques surfaces.

We will say that a divisor $D$ on a $K3$ surface satisfies the conditions 
$(**)$ if it satisfies 
the conditions $(*)$ and in addition $D^2 \geq 0$ and 
$(L^2,D^2) \not = (4k+2, k)$. 

Also recall that by a result of Green and Lazarsfeld \cite{gl} all smooth 
curves in a base point free linear system on a $K3$ surface have the same 
Clifford index (see Section \ref{cliff} for more details). The same is not 
true for the gonality (see Remark \ref{remgon} below).

\begin{thm} \label{birkvathm}
  Let $L$ be a globally generated, big line bundle on a $K3$ surface and 
  $k \geq 1$ an integer. Denote by $c$ the Clifford index of all smooth 
  curves in $|L|$.
  The following conditions are equivalent:
\begin{itemize}
  \item [a)] $L$ is birationally $k$-very ample,
  \item [b)] $L$ is birationally $k$-spanned,
  \item [c)] $L^2 \geq 4k$ and there exists no effective divisor $D$ 
  satisfying the conditions $(**)$, 
  \item [d)] $L^2 \geq 4k$ and there exists no smooth curve $D$ satisfying 
   the conditions $(**)$,
  \item [e)] the minimal gonality of a smooth curve in $|L|$ is $\geq k+2$,
  \item [f)] $c \geq k$,
   and there exists a smooth curve in $|L|$ having gonality $c+2$; or 
  $c=k-1$ and all smooth curves in $|L|$ have gonality $c+3$ (in which 
  case, $L \sim 2D + \Gamma$, where $D$ and $\Gamma$ are smooth 
  curves satisfying $D^2 = k$, $\Gamma ^2 =-2$ and $D.\Gamma=1$),
  \item [g)] there exists no smooth curve $C$ in $|L|$ containing a 
  $0$-dimensional subscheme of degree $k+1$ where the $k$-spannedness of 
  $L_{|C} \iso \omega_C$ fails,
  \item [h)] there exists no rank $2$ vector bundle $E$ on $S$ generated by 
  its global sections satisfying $H^1 (E)= H^2 (E) = 0$, $\det E = L$ and 
  $c_2 (E) \leq k+1$.
  \end{itemize}
\end{thm}

The paper is organized as follows.

In Section \ref{use} we fix notation and gather some results of Saint-Donat 
that will be needed in the
rest of the paper. Then we prove in Section \ref{numcond} that if 
$L^2 \geq 4k$ and $S$ is $K3$, or $L^2 \geq 4k+4$ and $S$ is Enriques, and 
$L+K_S$ is
not $k$-very ample (resp. not $k$-spanned), there exists a divisor $D$
containing some (resp. some curvilinear) 0-dimensional scheme $\Z$ of degree 
$\leq k+1$ where the $k$-very ampleness (resp. $k$-spannedness) fails and
satisfying certain numerical conditions (which are $(*)$ in the $K3$ case).

In Section \ref{failenr} we treat the Enriques case and give
the proof of Theorem \ref{kvathmenr}.

The rest of the paper will deal with $K3$ surfaces and the proofs of 
Theorems \ref{kvathm} and \ref{birkvathm}.

We first show, in Section \ref{smcur}, that it is always possible to find a 
smooth curve among all divisors satisfying $(*)$.

In Section \ref{fail} we conclude the proof of Theorem \ref{kvathm} by using such a smooth curve to show that $L$ is not 
$k$-spanned (and hence not $k$-very ample). We also investigate more closely the 
case $L^2 < 4k$. In any case we explicitly construct $0$-dimensional schemes of degree $k+1$ where the 
$k$-spannedness of $L$ fails. This explicit construction will be needed in Section \ref{birat}, where we show that the existence of any divisor satisfying $(**)$ will in fact imply that $L$ fails to be $k$-spanned on any Zariski-open subset of $S$, but if the only divisors $D$ satisfying $(*)$ are those satisfying the special conditions $D^2 \leq -2$ or 
$(L^2,D^2) = (4k+2, k)$, then $L$ is in fact birationally 
$k$-spanned, even though it is not $k$-spanned. This shows the equivalence of parts a)-d) in Theorem \ref{birkvathm}.

In Section \ref{cliff} we discuss the Clifford index and gonality of curves in $|L|$, relying upon results in \cite{gl} and \cite{cp}, and finish the proof of Theorem \ref{birkvathm}.

\begin{rem}
  {\rm Note that in \cite{dr2} Di Rocco showed that if $L$ is a $k$-very
  ample line bundle of degree $\leq 4k+4$ on a surface $S$ and $k \geq 2$, 
  then $(S,L)$ belongs to a certain
  list of pairs (\cite[Table 2]{dr2} \footnote{Note that in this list the case 
  of Enriques surfaces and line bundles of degrees $4k+4$ are  
  missing \cite{drpriv}. }), and all the 
  line bundles in the list are proved to be $k$-very ample, except
  for the $K3$ and Enriques cases. Thus, the results in this paper also 
  complete the description of $k$-very
  ample line bundle of degree $\leq 4k+4$ on a surface, for $k \geq 2$.}
\end{rem}

\begin{rem}
  {\rm Conditions for birational $k$-very ampleness and birational 
  $k$-spannedness can most probably be found for the other surfaces of 
  Kodaira dimension zero, but as far as we can see, the connections to 
  Clifford index and gonality of smooth curves in $|L|$ are not obtained as 
  easily as in the $K3$ case. In general, it would be interesting to know 
  whether there are connections between birational $k$-very ampleness and 
  birational $k$-spannedness and the Clifford index and gonality of smooth 
  curves in $|L|$ for other surfaces than $K3$s.}
\end{rem}

\begin{rem}
  {\rm Note that for $k=0$ and $1$ we retrieve the special results of 
  Saint-Donat concerning criteria for spannedness and very ampleness of 
  line bundles on $K3$ surfaces. Denoting by 
  $\Phi _{L}$ the morphism defined by the complete linear
  system $|L|$, when $L$ is spanned, we get the well-known result that 
  $\Phi _{L}$ is not birational if and only if all smooth curves in
  $|L|$ are hyperelliptic.}
\end{rem}

\begin{rem} \label{remgon}
  {\rm The only example known where the smooth curves in a base point free 
  linear system $|L|$ on a $K3$ surface do not have constant gonality, is the 
  famous Donagi-Morrison example \cite[(2.2)]{DM}. If $L$ is ample, Ciliberto 
  and Pareschi \cite{cp} have showed that this is indeed the only such 
  example, but the question remains open for the cases where $L$ is not ample.
  
  In addition to the Donagi-Morrison example, the only other example known of 
  exeptional curves in a base point free linear system on a $K3$ surface
  is an example of Eisenbud, Lange, Martens and Schreyer (see Remark 
  \ref{eisex}). This example appears in a natural way in our treatment of 
  birational $k$-very ampleness and birational $k$-spannedness (it is the 
  second case of part f) in Theorem \ref{birkvathm}), and in Section 
  \ref{cliff} we show that this is the only example of a base point free 
  linear system on a $K3$ surface where all smooth curves are exceptional.

  Unfortunately, we are not able, by our treatment, to ``explain'' the 
  Donagi-Morrison example in terms of birational $k$-very ampleness and 
  birational $k$-spannedness, nor to treat the question of the constancy of 
  gonality of the smooth curves in $|L|$ when $L$ is not ample. Those would 
  be very interesting questions to treat.}
\end{rem}

\noindent {\bf Acknowledgements.} I would like to thank S. Di Rocco, C. Hacon, P. Belkale, R. Piene 
for useful conversations, and H. Terakawa for interesting correspondence. Thanks are also due to the Department of Mathematics at the University of Utah for its hospitality and wonderful atmosphere during the time a considerable part of this work was done. 

Finally, I would like to thank the referee for informing me about the paper of Szemberg \cite{sz} and T.~Szemberg himself for useful comments during the writing of the final version of this paper.

\section{Notation and Background Material}
\label{use}

We use standard notation from Algebraic Geometry.

The ground field is the field of complex numbers. All surfaces are \textit{smooth algebraic surfaces}.

By a \textit{curve} on  a surface $S$ is always meant an \textit{irreducible curve} (possibly singular), i.e. a prime divisor. Line bundles and divisors are used with little or no distiction, as well as the multiplicative and additive notation. Linear equivalence of divisors is denoted by $\sim$, and numerical equivalence by $\eqv$. Note that on a $K3$ surface linear and numerical equivalence is the same.

If $L$ is any line bundle on a surface, $L$ is said to be \textit{numerically effective}, or simply \textit{nef}, if $L.C \geq 0$ for all curves $C$ on $S$. In this case $L$ is said to be \textit{big} if $L^2 >0$.

If $\F$ is any coherent sheaf on a variety $V$, we shall denote by 
$h^i (\F)$ the complex dimension of $H^i (V, \F)$, and by $\chi (\F)$ the Euler characteristic $\sum (-1)^i h^i (\F)$.

If $D$ is any divisor on a surface $S$, Riemann-Roch for $D$ is 
$\chi (\O_S (D))= \frac{1}{2}D.(D-K_S) + \chi (\O_S)$, where $K_S$ is the 
canonical bundle of $S$.

If $D$ is any effective divisor on $S$, and $\L$ any line bundle on $D$, 
Riemann-Roch yields $\chi (\L) = \deg \L + \chi (\O_D)=  \deg \L - 
\frac{1}{2}D(D+K_S)$.

By an Enriques surface is meant a surface $S$ with 
$H^1 (\O_S)=0$ and such that the canonical bundle $K_S$ satisfies $K_S \not 
\iso \O_S$, and $K_S ^2 \iso \O_S$. Recall that we also have $h^0 (K_S) 
= h^1 (K_S) = h^2 (\O_S)=0$, $h^2 (K_S) = h^0 (\O_S)= 1$ and $\chi (\O_S)= 1$.

By a $K3$ surface is meant a surface $S$ with trivial 
canonical bundle and such that $H^1 (\O_S)=0$. In particular $h^2 (\O_S)= 1$ 
and $\chi (\O_S)= 2$.

We will make use of the following results of Saint-Donat on line bundles on 
$K3$ surfaces. The first result will be used repeatedly, without further 
mention. 

\begin{prop}
\cite[Cor. 3.2]{S-D} A complete linear system on a $K3$ surface
has no base points outside of its fixed components.
\end{prop}

\begin{prop} \label{sd1}
  \cite[Prop. 2.6]{S-D} Let $L$ be an invertible sheaf on a $K3$
  surface $S$ such that $|L| \not = \emptyset$ and such that $|L|$ has no
  fixed components. Then either
  \begin{itemize}
  \item [i)] $L^2 >0$ and the generic member of $|L|$ is an irreducible
    curve of arithmetic genus $L^2/2 +1$. In this case $h^1 (L) =0$, or
  \item [ii)] $L^2=0$, then $L \iso \O _S (kE)$, where $k$ is an integer
    $\geq 1$ and $E$ is an irreducible curve of arithmetic genus 1. In
    this case $h^1 (L) = k-1$ and every member of $|L|$ can be written
    as a sum $E_1 + ... + E_k$, where $E_i \in |E|$ for $i=1, ... , k$. 
  \end{itemize}
\end{prop}

Note that by Bertini the generic members in $|L|$ and $|E|$ are smooth.

\begin{lemma} \label{sd2}
  \cite[2.7]{S-D} Let $L$ be a nef line bundle on a $K3$
  surface $S$. Then $|L|$ is not base point free if and only if there
  exist smooth irreducible curves $E$, 
    $\Gamma$ and an integer $k \geq 2$ such that 
    \[ L \sim kE+ \Gamma, \hspace{.2in} E^2=0, \hspace{.2in}  \Gamma^2=-2,
    \hspace{.2in} E.\Gamma=1. \] 
  In this case, every member of $|L|$ is of the form
  $E_1+...+E_k+\Gamma$, where $E_i \in |E|$ for all $i$. 
\end{lemma}

The following result was mentioned in the introduction. For $k$-very 
ampleness and $k \geq 2$, it is proved in \cite[Cor. 3.2]{ba}. The proof for 
the $k$-spanned case and $k \geq 0$ is not much more involved and therefore 
left to the reader.

\begin{prop} \label{under}
  Let $L$ be a big and nef line bundle on a smooth surface $S$ of Kodaira 
  dimension zero. If $L$ is $k$-spanned and $L^2 \leq 4k+4$, then either $S$ 
  is a $K3$ surface and $L^2= 4k$, $4k+2$ or $4k+4$, or $S$ is an Enriques 
  surface and $L^2= 4k+4$.
\end{prop}

\section{Numerical Conditions if $k$-very ampleness or $k$-spannedness
  is not Fulfilled}
\label{numcond}

We will now consider the case where $S$ is a $K3$ surface and $L^2 \geq 4k$, 
or $S$ is an Enriques surface and $L^2 \geq 4k+4$, and $K_S \* L$ fails to 
be $k$-very ample or $k$-spanned.

The following result is due to Beltrametti, Francia and Sommese.

\begin{prop} \label{bfs}
  Let $L$ be a nef and big line bundle on a surface $S$ and let $\Z$
be any $0$-dimensional subscheme of $S$ such that $\deg \Z = k+1$. Assume that the map
\[ H^0 (K_S \* L) \hpil H^0 (K_S \* L  \* \O _{\Z}) \]
is not onto, and for any proper subscheme $\Z'$ of $\Z$, the map 
\[ H^0 (K_S \* L) \hpil H^0 (K_S \* L  \* \O _{\Z'}) \]
is onto.  

Then there exists a rank 2 vector bundle $E$ on $S$ fitting
into the exact sequence 
\begin{equation} \label{BSes}
 0 \hpil \O_S \hpil E \hpil L \* \I _ {\Z} \hpil 0, 
\end{equation}
and such that the coboundary map of the exact sequence tensorized with $K_S$,
\[ \delta : H^1 (L \* K_S \* \I _ {\Z}) \hpil H^2 (K_S) \iso H^0 (\O_S) \iso \CC, \]
is an isomorphism.
\end{prop}

\begin{proof}
  This follows from the first part of the proof of
  \cite[Thm. 2.1]{BS2} and from \cite[(1.12)]{Ty}.
\end{proof}

\begin{cor}
  We have $h^1 (E \* K_S) \leq h^1 (\O_S)$ and $h^2 (E \* K_S)=0$.
  In particular, if $S$ is $K3$ or Enriques, we have $h^1 (E \* K_S)
  = h^2 (E \* K_S)=0$ (with $K_S \iso \O_S$ in the $K3$ case).
\end{cor}

\begin{proof}
  The first assertion is immediate from the short exact sequence above 
  (tensorized with $K_S$), 
  the second follows from the short exact sequence
\[ 0 \hpil L \* K_S \* \I _ {\Z} \hpil L \* K_S  \hpil L \* K_S \* \O _ {\Z} 
\hpil 0, \]
  and the fact that $h^1 (L \* K_S \* \O _ {\Z}) = h^2 (L \* K_S) =0$.
\end{proof}

Note that from the sequence (\ref{BSes}) we get $c_1 (E)^2 = L^2$ and 
$c_2 (E) = \deg \Z = k+1$. 

The approach in \cite{BFS} and \cite{BS2} is now based upon the fact that when 
$L^2 \geq 4k+5$, we have $c_1 (E) > 4 c_2 (E)$. By the well-known Bogomolov 
stability criterion (\cite{bog}, \cite{re}) one can then put $E$ in a 
suitable exact 
sequence. Since in our cases $L^2 \leq 4k+4$, we need a different approach.
Note that by our assumptions that $L^2 \geq 4k$ in the $K3$ case and 
$L^2 \geq 4k+4$ in the Enriques case, we have $c_1 (E)^2 - 4c_2 (E) 
\geq -4$ and $c_1 (E)^2 - 4c_2 (E) \geq 0$ in the $K3$ and Enriques case, 
respectively.

We will need the following result by Donagi and Morrison.

\begin{lemma} \label{domo}
  Let $E$ be a nonsimple, rank $2$ vector bundle on
  a surface $S$.
  There exist line bundles $M$, $N$ and a zero-dimensional subscheme $A
  \subset S$ such that $E$ fits in an exact sequence 
\[    0 \hpil N \hpil E \hpil M \* \I _ {A} \hpil 0 \]
  and either 
  \begin{itemize}
  \item [(a)] $N \geq M$, or
  \item [(b)] $A= \emptyset$ and the sequence splits.
  \end{itemize}
\end{lemma}

\begin{proof}
  This follows the lines of the proof of \cite[Lemma 4.4]{DM}, by noting 
  that the assumption that $S$ be a $K3$ is unnecessary.
\end{proof}

Now we prove the result which enables us to avoid the use of Bogomolov 
stability.

\begin{prop} \label{split}
  Let $E$ be a vector bundle of rank two on a $K3$ surface (resp. an 
  Enriques surface) satisfying $c_1 (E)^2 - 4c_2 (E) \geq -4$ (resp.
  $c_1 (E)^2 - 4c_2 (E) \geq 0$). Then there exist line bundles $M$, $N$ and a 
  zero-dimensional subscheme $A
  \subset S$ such that $E$ fits in an exact sequence 
  \begin{equation} \label{enes}
    0 \hpil N \hpil E \hpil M \* \I _ {A} \hpil 0
  \end{equation}
  and either 
  \begin{itemize}
  \item [(a)] $N \geq M$, or
  \item [(b)] $A= \emptyset$ and the sequence splits.
  \end{itemize}
\end{prop}

\begin{proof}
  For a vector bundle of rank $e$ on a surface $S$ Riemann-Roch
  gives
\[ \chi (E \* E^{*} ) = (e-1) ( c_1 (E))^2 - 2e c_2 (E) + e^2 \chi (\O_S). \]
  From our assumptions we get $\chi (E \* E^{*} ) \geq 4$ for both the $K3$ 
  and the Enriques case, whence  
  $h^0 (E \* E^{*} ) + h^0 (K_S \* E \* E^{*} ) \geq 4$. In the $K3$ case, 
  since $K_S$ is trivial, we get $h^0 (E \* E^{*} ) \geq 2$, and we are done 
  by Lemma \ref{domo}. So we will stick to the Enriques case.

  If $h^0 (E \* E^{*} ) \geq 2$, we are again done by Lemma \ref{domo}. 
  So we can 
  assume $h^0 (K_S \* E \* E^{*} ) \geq 3$ for the rest of the proof.

  Pick any non-zero section $\alpha: E \khpil K_S \* E$ in 
  $H^0 (K_S \* E \* E^{*}) \iso \Hom (E, K_S \* E)$. Define $\alpha' := 
  K_S \*  \alpha$ and consider the composite morphism
\[ \alpha' \circ \alpha : E \hpil E \* {K_S}^2 \iso E.  \]
  Then we have three possibilities:
\begin{itemize}
  \item [i)]   $\alpha' \circ \alpha$ is zero,
  \item [ii)]  $\alpha' \circ \alpha$ is a non-zero multiple of the identity,
  \item [iii)] $\alpha' \circ \alpha$ is not a multiple of the identity.
\end{itemize}

  In case iii), $\alpha' \circ \alpha$ gives a non-trivial endomorphism of 
  $E$, whence $E$ is non-simple and we are done again.

  In case ii) both $\alpha$ and $\alpha'$ have constant rank two, whence 
  $E \iso E \* K_S$ and $h^0 (E \* E^{*}) = h^0 (K_S \* E \* E^{*} ) \geq 3$, 
  and we are done again.

  In case i), $\alpha$ must drop rank. Since $\det \alpha \in \Hom (L,L) 
  \iso H^0 (L) \iso \CC$ and $\alpha$ is not zero, $\alpha$ has constant
  rank equal to one, whence $\ker \alpha =: N$ is a line bundle and we can 
  write $\im \alpha = M \* \I _ {A}$, where $M := (\im \alpha)^{**}$ 
  is a line bundle and $A$ is a zero-dimensional subscheme in of $S$. In 
  other words we have an exact sequence
\[    0 \hpil N \hpil E \hpil M \* \I _ {A} \hpil 0. \]

  Since $\alpha' \circ \alpha=0$, we have $\im \alpha \sub \ker \alpha \iso N
  \* K_S$, whence $ N \* K_S \* M \v$ has a section, and $N +K_S \geq M$. 
  From the short exact sequence above we get $c_1 (E) = M+N$ and $c_2 (E) = 
  M.N + \deg A$, whence 
\[ (M-N)^2 = c_1 (E)^2 - 4 M.N \geq c_1 (E)^2 - 4 c_2 (E) \geq 0. \]
  By Riemann-Roch either $N-M \geq 0$ or $M-N \geq 0$. If $M-N >0$, we 
  would get the absurdity 
  $K_S \geq M-N > 0$, whence $N-M \geq 0$, and we are done.
\end{proof}

\begin{rem} \label{enrrem}
 {\rm Actually, one can get a similar result in the case $c_1 (E)^2 - 4c_2 (E)
 = -2$ on an Enriques surface. In fact, the whole argument goes through, except if $h^0(E \* E^{*}) = 0$, $h^1(E \* E^{*})=1$ and $h^2(E \* E^{*}) =0$. Such a vector bundle is called} exceptional {\rm and by H. Kim \cite{kim} it sits in a
non-split short exact sequence
\[ 0 \hpil \O_S(D) \hpil E \hpil \O_S (D+\Gamma+K_S) \hpil 0, \]
for some divisors $D$ and $\Gamma$ such that $\Gamma ^2=-2$, 
$h^0(\Gamma)=1$ and $h^0(\Gamma+K_S)=0$. If $E$ is a vector bundle obtained as 
in Proposition \ref{bfs}, we have $D>0$.

We leave to the reader to work out the details. }
\end{rem}

Note that since $L$ is nef, we have $N.L \geq M.L$ in case (a), and in
case b) we can also assume this by symmetry. We will use this in the
proof of the next result, which is parallell to \cite[Prop. 1.4]{BFS}, but slightly different, due to the different hypotheses.

\begin{lemma} \label{numlem}
  With the same assumptions and notation as in Proposition \ref{bfs}, assume 
  furthermore that $S$ is a $K3$ surface and $L^2 \geq 4k$, or that $S$ is 
  an Enriques surface and $L^2 \geq 4k+4$.
  Let $E$ be the rank $2$
  vector bundle and $M$, $N$ the line bundles obtained as above and
  fitting in the sequences \rm{(\ref{BSes})} \textit{and}
  \rm{(\ref{enes})}. \textit{ Then the following
  conditions are satisfied (with $K_S=0$ in the $K3$ case):}
  \begin{itemize}
\item [i)]   \textit{The sequence $N \khpil E \khpil L \* \I _{\Z}$
             obtained from } \rm{(\ref{BSes})} \textit{and}
  \rm{(\ref{enes})} \textit{is nontrivial. 
\item [ii)]  $|M|$ contains an effective divisor D containing $\Z$.  
\item [iii)] $N >0$ and $H^1 (M+K_S) = H^2 (M) = H^2 (M+K_S) 
             = H^2 (N)=0$. Furthermore $M^2=D^2 \geq -2$.
\item [iv)]  $N - M \geq 0$ if $S$ is Enriques, or if $S$ is $K3$ and either
             $L^2 \geq 4k+2$ or $L^2= 4k$ and $A \not =
             \emptyset$.
\item [v)]   $(L-2D).L \geq 0$.
\item [vi)]  $L.D - \deg \Z \leq D^2$, with equality if and only if $A =
             \emptyset$.
\item [vii)] $D^2 \leq \frac{1}{2}L.D$, with equality if and only if $L
             \eqv 2D$.
\item [viii)]$D^2 \leq \deg \Z$, with equality if and only if $L
             \eqv 2D$ and  $L^2= 4k+4$.
\item [ix)]  $L.D \leq 2 \deg \Z$, with equality if and only if $L
             \eqv 2D$ and  $L^2= 4k+4$.
\item [x)]   If $(L.D)^2 = L^2 D^2$, then $L \eqv \lambda D$, for some
             $\lambda \in \QQ$, $2 \leq \lambda \leq 1 + \deg \Z / D^2.$}
  \end{itemize}
\end{lemma}

\begin{proof}
  \begin{itemize}
\item [i)] First we need to show that $h^0(-N)=0$. Note that $N$ is 
             non-trivial, since otherwise $N.L= M.L = L^2 =0$. Assume, to 
             get a contradiction, that $|-N|$ contains an effective member 
             $N_0$. Then, since $|L|$ contains an effective member (because 
             it is big and nef) and $L \sim M - N_0$, $|M|$ contains an 
             effective member $M_0$. But then $L.M_0 \leq L.N = - L.N_0
             \leq 0$ contradicts that $L$ is big and nef. 

             If the composition $N \khpil E \khpil L \* \I _{\Z}$ is zero,
             then $N \sub \ker (E \khpil L \* \I _{\Z})$, so $N=\O_S (-D')$ 
             for some $D'\geq 0$, contradicting the fact that $h^0(-N)=0$.
\item [ii)]  Tensorizing the sequences (\ref{BSes}) and (\ref{enes}) with
             $N^{-1}$ and taking cohomology, we get $0 \not = H^0 (E \* N^{-1})
             \sub H^0 (M \* \I _{\Z})$. 
\item [iii)] We first prove that $N >0$. This is clear if we are in 
             case a) of Lemma \ref{split}, since we have proved that $M>0$. 
             So we can assume that $A= \emptyset$ in sequence (\ref{enes}). 
             Since $L \sim M + N$ by (\ref{BSes}) and $N.L \geq M.L$, we
             find $N^2 \geq M^2$. From (\ref{BSes}) and (\ref{enes}) we 
             find $N.M = c_2 (E)= k+1$, so $N^2 \geq \frac{1}{2}(L^2-2M.N)
             \geq -1$ if $S$ is $K3$ and $N^2 \geq 0$ , if $S$ is Enriques. 
             In both cases, since $N$ is non-trivial, by Riemann-Roch either 
             $|N|$ or $|-N|$
             contains an effective member. So $N >0$, since we have 
             proved that $h^0(-N) =0$. By Serre duality 
             $h^2 (N) = h^0 (K_S-N) =0$.

  The exact sequence (\ref{split}) gives rise to an exact sequence
\[  0 \hpil N \* K_S \hpil E \* K_S \hpil M  \* K_S \hpil \tau \hpil 0, \]
 where $\tau$ is a torsion sheaf supported on a finite set. Taking
 cohomology and using that $H^1(E\* K_S ) = H^2(E\* K_S ) = 
 H^2 (N\* K_S ) = 0$, it follows that
 $H^1(M\* K_S ) = H^2(M \* K_S ) = 0$. Clearly we also have $h^2 (M) = 
 h^0 (K_S-M)=0$. By Riemann-Roch
\[ \chi (M \* K_S) = h^0 (M \* K_S) = \frac{1}{2} M^2 + \chi (\O_S), \]
whence $M^2=D^2 \geq -2$ in both the $K3$ and Enriques case.
\item [iv)]  Since this is clear when we are in case (a) of Lemma
             \ref{split}, we can assume we are in case (b) of that lemma,
             and hence that $\deg A=0$. From (\ref{enes}) we get
             $(N-M)^2 = L^2 - 4k - 4$. By our assumptions, we get
             $(N-M)^2 \geq 0$ if $S$ is Enriques and $(N-M)^2 \geq -2$ if 
             $S$ is $K3$, so by Riemann-Roch either
             $N-M \geq 0$ or $M-N \geq 0$. By
             symmetry, we can assume $N-M \geq 0$.
\item [v)]   This is $(N-M).L \geq 0$ rewritten.
\item [vi)]  By (\ref{BSes}) and (\ref{enes}) we have
             $c_2 (E) = (L-D).D + \deg A = \deg \Z$. 
\item [vii)] By Hodge Index Theorem $L^2 D^2 \leq (L.D)^2$, with equality
             if and only if $(D.L)L \eqv (L^2) D$, which means that $D
             \eqv aL$ for some $a \in \QQ^+$ (since $L^2 \not =0$ and $D
             >0$, one has $D.L \not =0$). 
             But since $L \eqv D + N$, this gives $N
             \eqv bL$ for some $b \in \QQ^+$. Combining this with v) we get 
             $2(D.L)D^2 \leq (L.D)^2$, with equality if and only if $D
             \eqv aL$, $N \eqv bL$ and $D.L = N.L$. This means that
             $a=b=1/2$ and $L \eqv 2D$.
\item [viii)]Immediate from vi) and vii).
\item [ix)]  Immediate from vi) and viii).
\item [x)]   This follows from v) and vi).
   \end{itemize}
\end{proof}

As a consequence, we have

\begin{prop} \label{numprop}
  With the same assumptions and notation as in Proposition \ref{bfs}, assume 
  furthermore that $S$ is a $K3$ surface and $L^2 \geq 4k$, or that $S$ is 
  an Enriques surface and $L^2 \geq 4k+4$.
  
  Then there exists an effective
  divisor $D$ passing through $\Z$ and such that
\begin{equation}
    \label{num1}
    L.D- k-1 \leq D^2 \stackrel{i)}{\leq} \frac{1}{2}L.D
    \stackrel{ii)}{\leq} k+1
\end{equation}
with equality in i) if and only if $L \eqv 2D$ and $L^2 \leq 4k+4$, and 
equality in ii) if
and only if $L \eqv 2D$ and $L^2= 4k+4$. 

Furthermore $L-2D \geq 0$ if $L^2 \geq 4k+2$.
\end{prop}

As a special case, we get that if $L+K_S$ is not $k$-very ample, then there 
exists an effective divisor $D$ as above.

The numerical conditions in (\ref{num1}) can also be formulated as
\begin{eqnarray*}
& 2D^2 \stackrel{i)}{\leq} L.D \leq D^2 + k+1 \stackrel{ii)}{\leq} 2k+2   \\
(*) & \mbox{ with equality in i) if and only if } L \sim 2D \mbox{ and } L^2 \leq 4k+4,   \\
& \mbox{and
equality in ii) if and only if } L \sim 2D \mbox{ and } L^2= 4k+4. 
     \end{eqnarray*}
if $S$ is a $K3$ surface, and as  \begin{eqnarray*}
(\#) & 2D^2 \stackrel{i)}{\leq} L.D \leq D^2 + k+1 \stackrel{ii)}{\leq} 2k+2   \\
& \mbox{ with equality in i) or ii) if and only if } L \eqv 2D \mbox{ and } L^2 = 4k+4. 
     \end{eqnarray*}
if $S$ is an Enriques surface.

\section{$k$-spannedness Failing on an Enriques Surface}
\label{failenr}

In this section we will study the Enriques case. Given a divisor satisfying 
$(\#)$, we will give an explicit construction of $0$-dimensional schemes 
where the $k$-spannedness of $L+K_S$ fails.

First of all, note (in both the Enriques and $K3$ cases) that if $D$ is an 
effective divisor satisfying the conditions $(\#)$ or
$(*)$ for $k=k_0$, and the middle inequality is strict, then $D$ will satisfy 
the same conditions for $k=k_0-1$. So if $D$ is a divisor satisfying the 
conditions $(\#)$ for the smallest integer $k_0$, then $D$ will have to satisfy
$D.L=D^2+k_0+1$.

Recall the following result due to Cossec and Dolgachev \cite[p. 134]{CD}:

\begin{thm} \label{cosdol}
  Let $D$ be an effective divisor with $D^2 >0$ on an Enriques surface. There 
  exists a divisor $f >0$ satisfying $f^2=0$ and $f.D \leq \sqrt{D^2}$. 
\end{thm}

As a consequence, we get:

\begin{lemma} \label{numred}
  Assume $D$ is an effective divisor satisfying the numerical 
  conditions $(\#)$ for an integer $k=k_0$ and that there are no divisors 
  satisfying the conditions for any integer $k < k_0$.

  Then $D$ is of one of the following types:
\begin{itemize}
\item[i)]   $D^2=-2$, $D.L=k_0-1$,
\item[ii)]  $D^2=0$,  $D.L=k_0+1$,
\item[iii)] $D^2=2$, $L \eqv 2D$, $k_0=1$,
\item[iv)]  $D^2=4$, $L \eqv 2D$, $k_0=3$,
\end{itemize}
and if $D$ is as in iii) or iv), then $D \sim E_1+E_2$, where $E_1^2=E_2^2=0$ 
and $E_1.L=E_2.L=k_0+1$.
\end{lemma}

\begin{proof}
  We clearly must have $D^2 \geq -2$. Also note that the result is clear for 
  $k_0 \leq 1$, so we will assume $k_0 \geq 2$.

  The divisor $D$ has to be nef if $D^2 \geq 0$. Indeed, if there existed a curve $\Gamma$ 
  such that $\Gamma ^2=-2$ and $\Gamma.D <0$, then $\Gamma.L=0$ (see the proof 
  of Claim 2 in Proposition \ref{smcurprop} below) and $\Gamma$ would satisfy 
  the conditions $(\#)$ for $k=1$, a contradiction.

  If $D ^2 \geq 2$, then by Theorem \ref{cosdol} we can find an 
  $f >0$ such that $f^2=0$ and such that 
  $(D-f)^2 \geq D^2- 2\lfloor \sqrt{D^2} \rfloor \geq 0$ and 
  $(D-f).D \geq D^2- \lfloor \sqrt{D^2} \rfloor > 0$. By Riemann-Roch and the 
  nefness of $D$, we have that $A := D-f >0$. Since 
  $L.D = L.f +L.A \leq 2k_0+2$, 
  we must have $L.f=L.A=k_0+1$, $A^2=0$, $D^2=2$ or $4$ and 
  $L.D= D^2 + k_0 + 1=2k_0+2$, 
  whence $L \eqv 2D$ and $k_0=1$ or $3$ respectively.
\end{proof}

We will now prove

\begin{prop} \label{enrfail}
  Let $L$ be a big and nef line bundle on an Enriques surface $S$ satisfying 
$L^2 \geq 4k_0+4$. Assume there exists a divisor $D$ satisfying the numerical 
conditions $(\#)$ for an integer $k=k_0$ and that there are no divisors 
satisfying the conditions for any integer $k<k_0$ (so that, in particular, 
$D$ must be as in i)-iv) in Lemma \ref{numred} above).

Define $F := L-D$ and $F_D := F \* \O_D$. Then the generic section of $F_D$ 
will define a zero dimensional curvilinear scheme of degree $k_0+1$ where the 
$k_0$-spannedness of $(L+K_S)_D := (L+K_S) \* \O_D$ fails. 
\end{prop}

\begin{proof}
  The numerical conditions above give $D^2 \leq k_0+1$, with equality if and 
  only if $L \eqv 2D$. Note also that $F.D= k_0+1$

  Since $D^2 \leq k_0+1$, we get by Riemann-Roch on $D$,
\[ h^0 (F_D) \geq k_0+1 - \frac{1}{2} D^2 \geq \frac{1}{2}k_0 + \frac{1}{2} >0.\]
  Any non-zero section $s$ in $H^0 (F_D)$ gives a short exact sequence (after 
  tensorising with the dualising sheaf $\omega ^{\circ}_D := K_S \* \O_D (D)$)
\[ 0 \khpil \omega ^{\circ}_D  \khpil (L+K_S)_D \ \khpil G \khpil 0,\]
  where $G$ is a torsion sheaf supported on the zero-scheme of $s$, so 
  $\length G = \deg F_D = k_0+1$.

  Now $h^1 ((L+K_S)_D) = h^0 (-F_D) = 0$ and $h^1 (\omega ^{\circ}_D) = 
  h^0 (\O_D) >0$, so $(L+K_S)_D$ fails to be $k_0$-very ample on the zero 
  scheme of $s$. To show that $(L+K_S)_D$ is not even $k_0$-spanned, we can 
  assume $k_0 \geq 2$. 

  If $D$ is as in case iv) of Lemma \ref{numred}, then 
  $F \eqv D \eqv \frac{1}{2}L$ is nef. If $|F|$ is not base point free, 
  then since $F^2 =D^2 =4$, there exists by Proposition \ref{numprop} an 
  effective 
  divisor $B$ satisfying $B^2=0$ and $B.F=1$. But this would imply $B.L=2$, 
  and $B$ would satisfy the conditions $(\#)$ for $k=1 < k_0$, a 
  contradiction. So $|F|$ is base point free and the zero scheme of a 
  generic section of $F_D$ 
  is curvilinear. Hence we are done.
  
  By Lemma \ref{numred}, what remains are the cases where $D$ is as in i) 
  or ii). In particular, $L.D \leq k_0+1$.

  We will show that $D$ is reduced and that $F_D$ is base point 
  free on $D$. It then follows that the zero scheme of a generic section of 
  $F_D$ is curvilinear, which will complete the proof.

  If $D$ is reduced and irreducible, then since 
  $\deg F_D = k_0+1 \geq 2 p_a(D) = D^2+2$, we have that $F_D$ is base point 
  free on $D$ (see e.g. \cite[Prop. 2.3]{CF}).

  So assume $D = \sum_{i=1}^n D_i$, for an integer $n \geq 2$, with all 
  $D_i >0$. Then $L.D = \sum_{i=1}^n L.D_i \leq k_0+1$, and since
  none of the $D_i$ can satisfy the conditions $(\#)$ for any $k<k_0$ by 
  assumption, we easily find $k_0=2$ or $3$, $n \leq 3$, $D^2=0$, $L.D=k_0+1$,
  $D_i^2 =-2$ and 
  $L.D_i \geq 1$ for all $i$. Since $D^2=0$, we have that $D$ is nef, as in 
  the proof of Lemma \ref{numred}, whence $D.D_i=0$ for all $i$. This gives  
  $\deg F_{D_i} = (L-D).D_i \geq 1 > 2 p_a(D_i)$. By 
  \cite[Prop. 2.3]{CF} again, $F_D$ is base point free on $D$.

  Since $D_i^2=-2$ for all $i$ (and any decomposition), $D$ is reduced.
\end{proof}

 \begin{rem}
  {\rm Note that one can argue in the same way (i.e. show that $L+K_S$ fails 
  to be $k$-very ample on any section of $F_D$ and show that a generic such 
  section is curvilinear) in the case $L^2 \geq 
  4k+5$ for surfaces of Kodaira dimension zero in general, and get a 
  simpler proof of \cite[Thm. 5]{ter}. Actually, 
  this approach gives the stronger result that if  $L$ is a big and nef 
  divisor on any surface $S$ satisfying $L^2 \geq 4k+5$ and 
  there exists an effective divisor $D$ satisfying $L.D - k-1 \leq D^2 < 
  L.D/2 < k+1$ and $D.K_S \leq 0$, then $K_S+L$ is not $k$-spanned 
  (here we leave some verifications to the reader). So Del Pezzo surfaces
  will for instance also be included.
  Note that $k$-very ample line bundles on Del Pezzo surfaces have been 
  completely characterized in \cite{dr1}.}
\end{rem}

Since $L$ and $L+K_S$ are numerically equivalent, we have shown for 
$L^2 \geq 4k+4$, that $L$ is $k$-very ample if and only if it is $k$-spanned if 
and only if there are no divisors satisfying the conditions $(\#)$. And among these divisors, we can always find one satisfying the conditions i) and ii) in Lemma \ref{numred}, by the last statement in that lemma.

To finish the proof of Theorem \ref{kvathmenr}, it suffices to note that if 
$L^2 \leq 4k+2$, there exists by Theorem \ref{cosdol} an $f >0$ 
satisfying $f^2=0$ and $f.L \leq \lfloor \sqrt{4k+2} \rfloor \leq k+1$. 

Note that we have also proved the following:

\begin{cor} \label{enrcor}
  If $L$ is a line bundle on an Enriques surface satisfying $L^2 \geq 4k+4$ and $L$ is $(k-1)$-very ample but not $k$-very ample, then any $0$-dimensional scheme $\Z$ where the $k$-very ampleness of $L$ fails is contained in an effective divisor $D$ satisfying one of the conditions i)-iv) in Lemma \ref{numred},
and any such divisor $D$ will contain $0$-dimensional schemes where the 
$k$-spannedness of $L$ fails.
\end{cor}

As a consequence of Theorem \ref{cosdol}, a line bundle of degree $4k+4$ cannot be $k$-very ample unless $k=0$, as remarked also by T. Szemberg in \cite{sz}.

However, our approach shows that also for $L^2=4k+4$ and $k \geq 1$, any (minimal, in the sense of Proposition \ref{bfs}) $0$-dimensional scheme $\Z$ where the $k$-very ampleness of $L$ fails, is contained in a divisor $D$ satisfying the conditions above. (By Remark \ref{enrrem} a similar result holds for $L^2=4k+2$.)

As an easy example, a line bundle of degree $8$ is not very ample (as was 
already known to Okonek \cite{ok}), and Corollary \ref{enrcor} shows that (if $L$ is base point free) each length $2$ scheme $\Z$ that $L$ fails to separate actually lies on a divisor $D>0$ satisfying i) $D^2=-2$, $D.L=0$, or ii) $D^2=0$, $L.D=2$, or iii) $D^2=2$, 
$L \eqv 2D$. Conversely, given such a divisor, Proposition \ref{enrfail} shows how to explicitly find $0$-dimensional schemes of length $2$ that $L$ does not separate. By Theorem \ref{cosdol} a divisor of type ii) above always exists. The existence of a divisor of type i) corresponds to $L$ contracting rational curves, whereas the existence of a divisor of the type iii) is crucial to determine whether $\phi_L$ is birational or not (see \cite[Lemma 5.2.7]{cos}).

The next step in the paper is to prove the corresponding result to Proposition \ref{enrfail} when $S$ is a 
$K3$ surface and $L^2 \geq 4k$. Because of the weaker conditions, we will use a different approach and start by finding \textit{smooth 
curves} satisfying the conditions $(*)$.

\section{Smooth Curves Satisfying the Conditions $(*)$}
\label{smcur}

For the rest of the paper, we will consider the case where $S$ is a $K3$ 
surface.

We will now show that we can always find a smooth irreducible curve among
the divisors satisfying $(*)$. We will actually prove a stronger
result, which we will need in the next sections.

Recall that we commented in the beginning of the previous section that 
if $k_0$ is the smallest integer such that there exists a divisor $D$ 
satisfying $(*)$, then $D$ will have to satisfy $L.D=D^2+ k_0 +1$.

In exactly the same way it follows that if $D$ is any divisor
with $D^2 \geq 0$ satisfying the numerical conditions $(*)$ for
$k=k_0$ and there are no divisors with non-negative self-intersection
satisfying the numerical conditions $(*)$ for $k < k_0$, then $D$ will 
have to satisfy $L.D=D^2+k_0+1$.

With this in mind we prove

\begin{prop} \label{smcurprop}
  Let $S$ be a $K3$ surface and $L$ a big and nef line bundle on $S$, with
  $L^2 \geq 4k_0$ for an integer $k_0 \geq 0$. Assume there is a divisor
  $D$ satisfying the numerical conditions $(*)$ for 
  $k=k_0$. 
\begin{itemize}
  \item [a)] If $D^2 \leq -2$, there exists a smooth rational curve
    $D_0$ satisfying the conditions $(*)$ for $k=k_0$.
  \item [b)] If $D^2 \geq 0$ and furthermore there exists no divisor
    with non-negative self-intersection satisfying $(*)$
    for $k < k_0$ (so that, in particular, $L.D= D^2+k_0+1$),
    there exists a smooth curve 
    $D_0$ satisfying $0 \leq D_0^2 \leq D^2$, $L.D_0= D_0^2+k_0+1$ and the
    conditions $(*)$ for $k=k_0$, and furthermore  
    such that $F := L-D_0$ satisfies $h^0(F) \geq h^0(D_0)$.
\end{itemize}
\end{prop}

\noindent \textit{Proof.}
  If $D^2 \leq -2$, $D$ must have at least one smooth rational curve $D_0$
  as its component. Since $L$ is nef, this curve will clearly satisfy
  $L.D_0 \leq L.D$, so we get
\[ -4 = 2 {D_0}^2 < L.D_0 \leq L.D \leq D^2 +k+1 \leq  {D_0} ^2
  +k+1 = k-1, \]
  so $D_0$ is the desired curve. This proves a).

In case b), we first show with the help of three claims that we can
reduce to the case where $L^2 \geq 2D.L$, $H^1(D) =0$ and $D$ is nef
\footnote{Note that the two first properties are fulfilled if $D$ is a
  divisor obtained as in Lemma \ref{numlem}.}.

\vspace{.4cm}

\textbf{Claim 1.}
  \textit{We can assume $L^2 \geq 2D.L$ (equivalently $F^2 \geq D^2$ or
    $F.L \geq D.L$).}

\renewcommand{\proofname}{Proof of claim}

\begin{proof}
 From $(*)$ we have $2D.L \leq 4k_0 +4$. If equality holds,
 then $L \sim 2D$, so we have $L^2 =2D.L$. If $2D.L \leq
 4k_0$, or $L^2 \geq 4k_0 +2$, we also have  $L^2 \geq 2D.L$. So we only
 have the case $2D.L = 4k_0 +2$ and $L^2 =4k_0$ left. Since $L$ is big, we must have $k_0 \geq 1$, and from $L.D = D^2 + k_0 + 1$, we get $D^2 =k_0$, whence $k_0 \geq 2$. We have $L^2 -
 2D.L = -2$, which is equivalent to $D.L = F.L +2$ or $L^2 = 2F.L + 2$,
 for $F := L-D$. Since $L \sim D+F$ we also get $D^2=F^2+2$. In
 particular $F^2 \geq 0$. By Hodge Index $F^2 L^2 \leq (F.L)^2$, whence
 $2F^2 < F.L$. Furthermore
\[ F.L = D.L -2 = D^2+k_0-1 = F^2+k_0+1  \leq 2k_0, \]
 so we can interchange $D$ with $F$.
\end{proof} 

\textbf{Claim 2.}
  \textit{We can assume $D$ is nef.}
  \begin{proof}
If $\Gamma$ is a smooth rational curve (necessarily contained in the
base locus of $|D|$) such that $\Gamma.D <0$, define $D':= D-\Gamma$.
Then we have $D'.L \leq D.L \leq  \frac{1}{2} L^2$, whence by Hodge
Index $2{D'}^2 \leq L.D'$, with equality if and only if $L \sim 2D'$.

Furthermore we have
\[ {D'}^2 = D^2 - 2D.\Gamma + \Gamma ^2 =  D^2 - 2D.\Gamma -2 \geq D^2,\]
whence 
\[ D'.L \leq D.L = D^2+ k_0 +1 \leq  {D'}^2+ k_0 +1. \]
From the last equation combined with $2{D'}^2 \leq L.D'$, we get ${D'}^2
\leq k_0 +1$, with equality if and only if $L \sim 2D'$.

By the assumption that there exists no divisor
    with non-negative self-intersection satisfying $(*)$
    for $k < k_0$, we see that we must have $D'.L = D.L$ and ${D'}^2 = D^2$.
Combining all this, we see that we can exchange $D$ with $D'$, and since
$|D'|$ has one base divisor less than $|D|$, we are done by induction on
the number of base components of $|D|$, counted with multiplicities.
  \end{proof}

\textbf{Claim 3.}
  \textit{We can assume that $h^1 (D)=0$.}

  \begin{proof}
    By choosing $D$ as in Claim 2, it follows by Proposition \ref{sd1}
    and Lemma \ref{sd2} that $h^1 (D)=0$, unless when $D \sim lE$,
    for an integer $l \geq 2$ and $E$ a smooth elliptic curve. But in
    this case $L.E = \frac{1}{l} L.D = \frac{1}{l} (k_0+1)$, so $E$
    would satisfy the conditions $(*)$ for some $k < k_0$,
    contrary to our assumptions. 
  \end{proof}

\renewcommand{\proofname}{Rest of proof of Proposition 
              {\rm \ref{smcurprop}b)}}

\begin{proof}
By Proposition \ref {sd1} and Lemma \ref{sd2}, choosing $D$ according to the
claims above, the generic member of $|D|$ is a smooth curve of genus
$\geq 1$, unless 
\[ D \sim lE+ \Gamma,\] where $E$ and $\Gamma$ are smooth irreducible
curves satisfying $E^2=0$, $\Gamma^2=-2$ and $E.\Gamma=1$ and $l \geq 2$ an integer.

By the conditions $(*)$, we have $L.D \leq 2k_0+2$, so we get
\[ E.L \leq \frac{1}{l}D.L \leq \frac {2(k_0+1)}{l} \leq k_0+1, \]
and $E$ will satisfy the conditions $(*)$ for some $k \leq k_0$,
so we either have a contradiction, or $E$ is the desired smooth curve.

We now choose a smooth curve $D_0$ as above and we must show that $F :=
L-D_0$ has the desired property. 

Since $F^2 \geq D_0 ^2 \geq 0$, by Riemann-Roch either $F \geq 0$ or 
$-F \geq 0$. Since $F.L \geq D_0.L = k_0+1 >0$, we have $F >0$. By Riemann-Roch again, using $h^1(D_0) =0$:
\[ h^0(F) = \frac{1}{2} F^2 +2 + h^1 (F) \geq \frac{1}{2} D^2 +2 = h^0 (D). \]

This concludes the proof of Proposition \ref{smcurprop}
\end{proof}

\renewcommand{\proofname}{Proof}

\section{$k$-spannedness Failing on $K3$ Surfaces}
\label{fail}

We now conclude the proof of Theorem \ref{kvathm} by explicitly constructing, in the next two propositions, $0$-dimensional curvilinear schemes where the $k$-spannedness of $L$ fails, given the assumptions that $L^2 < 4k$ or that $L^2 \geq 4k$ and the existence of a divisor satisfying the conditions $(*)$ above.

We will make use of the following result by Beltrametti and Sommese:

\begin{thm} \label{bs12}
  \cite[(1.2)]{BS1} If $L$ is a special $k$-spanned line bundle on
  a smooth curve of genus $g$, then $g \geq 2k+1$.
\end{thm}

We first study the case $L^2 < 4k$ more thoroughly.

\begin{prop} \label{minus}
  Let $L$ be a globally generated, big line bundle on a $K3$ surface $S$
  satisfying $L^2 < 4k$ for an integer $k \geq 1$.
  Then any smooth curve $C \in |L|$ will contain a base point free complete 
  linear system $|A|$ of dimension $\geq 1$ and of degree $l+1$, for an integer  $1 \leq l \leq k$, such that the 
  $l$-very ampleness of $L$ fails on each member of $|A|$.
\end{prop}

\begin{proof}
  Let $C$ be any smooth curve in $|L|$. Then $g(C) = \frac{1}{2}L^2 +1
  \leq 2k$, so by Theorem \ref{bs12}, $L_C := L \* \O_C = \omega _C$
  is not $k$-spanned on $C$. So there exists a $0$-dimensional subscheme 
  $\Z$ of $C$ of degree $l+1$, for some integer  $1 \leq l \leq k$, such that the map
\[ H^0 (\omega _C) \hpil  H^0 (\omega _C \* \O_{\Z}) \]
  arising from the short exact sequence
\[ 0 \hpil \omega _C \* \O_C (-\Z) \hpil \omega _C \hpil \omega _C \*
\O_{\Z} \hpil 0 \]
  is not surjective. 

  Taking cohomology and using Serre duality, this is equivalent to
\[ h^1 (\omega _C \* \O_C (-\Z)) = h^0 (\O_C (\Z)) \geq 2. \]
 
  It is also clear that we can pick a minimal such $\Z$, i.e. with the 
  property  that the map  
  \[ H^0 (\omega _C) \hpil  H^0 (\omega _C \* \O_{\Z'}) \] 
  is surjective for all proper subschemes $\Z'$ in $\Z$. Then $|\O_C (\Z)|$ 
  will be base point free, for if not, we would by removing base points, 
  get a contradiction to the minimality of $\Z$.
\end{proof}

\begin{prop} \label{notconrat}
 With the same assumptions and notation as in Proposition \ref{smcurprop}, define $F:= L-D_0$. Then, for all smooth curves $D'_0$ in $|D_0|$ the line bundle $F_{D'_0} := F \* \O_{D'_0}$ has the property that the $k_0$-spannedness of $L_{D'_0}:= L\* \O_{D'_0}$ fails on any member of $|F_{D'_0}|$. 
\end{prop}

\begin{proof}
  This follows the first lines of the proof of Proposition \ref{enrfail}. 
  Since $D_0$ is smooth, any member of $|F_{D_0}|$ will automatically 
  be a curvilinear scheme where the $k_0$-spannedness of $L$ fails.
\end{proof}

This concludes the proof of Theorem \ref{kvathm}.

In the sequel we will need the following observation.

\begin{lemma} \label{smcurlem}
  With the same notation and assumptions as in Proposition \ref{notconrat}, assume furthermore that $L$ is generated by its global sections (so that in particular $k_0 \geq 1$) and that $D^2 \geq 0$ (and hence ${D_0}^2 \geq 0$). 

Then $h^0 (F_{D'_0}) \geq 2$ for any smooth curve $D'_0$ in $|D_0|$.

If furthermore $(L^2,D^2) \not = (4k_0+2, k_0)$, then $|F_{D'_0}|$ is base point free for any smooth curve $D'_0$ in $|D_0|$.
\end{lemma}

\begin{proof}
  First note that if $(L^2,D^2) \not = (4k_0+2, k_0)$, then the smooth curve $D_0$ constructed as in Proposition \ref{smcurprop} will also satisfy $(L^2,{D_0}^2) \not = (4k_0+2, k_0)$.

To alleviate notation, we will work with $D_0$. We first prove the statement that $h^0 (F_{D_0}) \geq 2$.

If $L \sim 2D_0$ and $L^2=4k_0+4$, then $F_{D_0} \iso \omega _{D_0}$ and $h^0 (\omega_{D_0}) = \frac{1}{2} k_0 + \frac{3}{2} \geq 2$.

In the other cases, since we can assume that 
${D_0}^2 \leq k_0 $, we get by Riemann-Roch that 
$h^0 (F_{D_0}) \geq  k_0 + 1 - \frac{1}{2} {D_0}^2 \geq   k_0 + 1 - 
  \frac{1}{2} k_0  > 1$.

Next we treat the question of the base point freeness of $|F_{D_0}|$.

\textbf{Case I:} $L \sim 2D_0$.

In this case, since $g(D_0) \geq 2$, we clearly have that $F_{D_0} \iso \omega _{D_0}$ is base point free.

\textbf{Case II:} $L^2 =4k_0$, $L \not \sim 2D_0$.

Since $L^2 = F^2 + D_0 ^2 + 2k_0 +2$ and $F^2 \geq D_0 ^2$, we must have
$D_0 ^2 \leq k_0-1$. Then
\[ \deg F_{D_0} = k_0+1 \geq D_0 ^2 +2 = 2g(D_0),\]
so $|F_{D_0}|$ is base point free \footnote{In fact, one can also prove that $|F|$ is base point free on $S$ in this case.}.

\textbf{Case III:} $L^2 \geq 4k_0+2$, $L \not \sim 2D_0$. 

In this case we have 
\[ (F-D_0)^2 = L^2 - 4k_0-4 \geq -2,\]
so by Riemann-Roch again, either $|F-D_0|$ or $|D_0-F|$ contains an
effective member (recall that we are assuming $F \not \sim D_0$). 

If $F.L>D_0.L$, then $|F-D_0|$ contains an
effective member. If $F.L=D_0.L$, then since $F \not \sim D_0$, and the argument for finding the smooth curve $D_0$ in Proposition \ref{smcurprop}
is symmetric, we can assume that it is $|F-D_0|$ that contains an
effective member.

Letting $\Delta$ be the (possibly zero) base divisor of $|F|$, we can
write
\[ F \sim D_0 + D_1 + \Delta,\]
for some divisor $D_1 \geq 0$.

 If either $D_0.D_1 >0$ or $\Delta.D_0 \geq 2$, we have 
\[ \deg \O_X (F) \* \O_{D_0} \geq D_0 ^2 +2 = 2g(D_0), \]
so $\O_X (F)$ is base point free on $D_0$, and we are done.

We now investigate the case $D_0.D_1 =0$ and $\Delta.D_0 =1$, separating it into the two cases $D_1=0$ and $D_1>0$.

\textbf{Case IIIa):} $D_1=0$.

This gives $F \sim D_0 + \Delta$ and $L \sim 2D_0 + \Delta$. Since $L$
is nef, we must have $L.\Delta = 2D_0. \Delta + \Delta ^2 \geq 0$, whence
$\Delta ^2 = -2$. 

Since $D_0$ is irreducible of non-negative self-intersection, the support of 
$\Delta$ can only meet $D_0$ once, so there has to exist a smooth rational 
curve $\Gamma$, occurring as a component of $\Delta$ of multiplicity one, such that $\Gamma.D =1$.

Again, since $L$ is nef, we must have $L.(\Delta - \Gamma) =
\Delta.(\Delta - \Gamma) = \Delta ^2 - \Delta. \Gamma \geq 0$. This
gives $\Delta. \Gamma \leq \Delta ^2 = -2$. We then calculate  
\[ (\Delta - \Gamma)^2 = -4 - 2 \Delta.\Gamma \geq 0. \]

If $\Delta \not = \Gamma$, then $(\Delta - \Gamma)$ would be contained
in the base locus of $|F|$, whence the contradiction $(\Delta -
\Gamma)^2 <0$.

So $\Delta = \Gamma$, and we calculate $({L_0}^2,{D_0}^2) = (4k_0+2, k_0)$. Since we easily see that if $(L^2,D^2) \not = (4k_0+2, k_0)$, then the smooth curve $D_0$ constructed as in Proposition \ref{smcurprop} will also satisfy $({L_0}^2,{D_0}^2) \not = (4k_0+2, k_0)$, we get that we are in the particular case 
$(L^2,D^2) = (4k_0+2, k_0)$ above.

\textbf{Case IIIb):} $D_1 > 0$.

Since $D_0.D_1 =0$ and $|D_0|$ is base point free, we can pick two
effective divisors respectively in $|D_0|$ and $|D_1|$ which do not
meet. Their sum is then an effective divisor in $|F-\Delta|$ which is
not 1-connected, whence $h^1(F-\Delta) >0$.

By Proposition \ref {sd1}, this means that $F-\Delta \sim lE$, for an
integer $l \geq 2$ and $E$ a smooth elliptic curve, and all members of
$|F-\Delta|$ would be a sum of elliptic curves in $|E|$. Hence $D_0 \sim
E$ and we can write $L \sim (l+1)E + \Delta$, with $E.\Delta = D_0.\Delta = 1$. But then $L.E=1$, and this contradicts the spannedness of $L$ by Theorem \ref{kvathm}.
\end{proof}

\section{Birational $k$-very ampleness and Birational $k$-spannedness}
\label{birat}

We have now found criteria for a big and nef line bundle $L$ on a 
$K3$ surface or Enriques surface $S$ to be $k$-very ample or equivalently 
$k$-spanned. But we still have not addressed the question about ``how much'' 
$L$ fails to be $k$-very ample or $k$-spanned. Is it possible to find a 
Zariski-open subset of $S$ where $L$ is $k$-very ample (resp. $k$-spanned), 
or does the $k$-very ampleness (resp. $k$-spannedness) of $L$ fail in 
such a ``bad'' way that such an open set is impossible to find ?

These questions motivate our further discussion and the following definition:

\begin{defn}  \label{birdef}
  Let $L$ be a globally generated line bundle on a smooth connected surface 
  $S$ and $k \geq 1$ an integer. 

  $L$ is \textit{birationally $k$-very ample} (resp. \textit{birationally 
  $k$-spanned)}, if there exists a non-empty Zariski-open 
  subset of $S$ where $L$ is $k$-very ample (resp. $k$-spanned).
\end{defn}

Recall that if $L$ is $0$-very ample (i.e. generated by global
sections), then the complete linear system $|L|$ defines a morphism $S
\khpil \PP ^{h^0 (L) -1}$. If $L$ is ($1$-)very ample, this morphism is
an embedding, and if $L$ is birationally ($1$-)very ample, this morphism is
birational. 

We will say that a divisor $D$ on a $K3$ surface satisfies the conditions $(**)$ for some integer $k \geq 1$, if $D$ satisfies $D^2 \geq 0$, $(L^2,D^2) \not = (4k+2, k)$ and the conditions $(*)$ above. From the results in the previous paragraph, we get the following:

\begin{prop} \label{firstbir}
  Let $L$ be a spanned and big line bundle on a $K3$ surface. If $L^2 < 4k$ or 
  there exists an effective divisor $D$ satisfying the conditions $(**)$, 
  $L$ is not birationally $k$-spanned.
\end{prop}

\begin{proof}
  We showed in the previous paragraph that given any of the assumptions above 
  there exists a complete linear system on $S$ (being either $|L|$ or 
  $|D_0|$) of dimension $\geq 1$, where the generic member is a smooth 
  curve each containing a base 
  point free linear system of dimension $\geq 1$ such that $L$ fails to 
  be $k$-spanned on each member of this linear system. Since any 
  Zariski-closed proper subset of $S$ will contain at most finitely many 
  of these curves and 
  intersect the rest of them in a finite number of points, the assertion  
  follows.
\end{proof}

What we have to consider now, is what happens in the case where there exist 
divisors $D$ satisfying $(*)$ and $D^2 < 0$ or $(L^2,D^2) = (4k+2, k)$.

\begin{lemma} \label{neg}
  Let $L$ be a globally generated, big line bundle on a $K3$ surface 
  $S$ satisfying $L^2 \geq 4k_0$ for an integer $k_0 \geq 0$. If the only 
  divisors satisfying $(*)$ for $k=k_0$ have negative 
  self-intersection, then $L$ is birationally $k_0$-very ample. 
\end{lemma}

\begin{proof}
  Assume that $\Z$ is any $0$-dimensional subscheme of $S$ of length $l+1$, 
  for an integer $l \leq k_0$, such that the map 
\[ H^0 (L) \hpil H^0 (L  \* \O _{\Z}) \]
  is not onto, but for any proper subscheme $\Z'$ of $\Z$, the map 
\[ H^0 (L) \hpil H^0 (L  \* \O _{\Z'}) \]
  is onto, then by the results in Section \ref{numcond}, there exists an 
  effective divisor $D$ containing $\Z$ and satisfying the numerical 
  conditions $(*)$ for $k=l$, and also $D^2 \geq -2$ and $h^1 (D)=0$. By our 
  assumptions,
  the divisor $D$ must satisfy $D^2=-2$ and consequently by Riemann-Roch, 
  $h^0 (D)=1$, whence $D$ is supported on a union of smooth rational curves. 
  In addition any such rational curve $\Gamma$ must satisfy 
  $\Gamma.L \leq D.L \leq l-1 \leq k_0 -1$.

  This means that any $0$-dimensional subscheme of $S$ where the $k_0$-very 
  ampleness of $L$ fails, must have some points in common with a smooth 
  rational 
  curve of degree $\leq k_0+1$ (with respect to $L$). Since the number of such 
  curves is finite by standard arguments, $L$ will be $k_0$-very ample on 
  the complement of the set of such curves, which is a Zariski-open subset 
  of $S$. 
\end{proof}

The case when $(L^2,D^2) = (4k+2, k)$ is a bit more involved. By the 
conditions $(*)$ we see that such a divisor $D$ must satisfy $D.L = 2k+1$.

\begin{lemma} \label{help1}
  Let $D$ and $L$ be effective divisors on a $K3$ surface satisfying 
  $L-2D \geq 0$, $L^2 =4k+2$, $D^2=k$ and  $D.L = 2k+1$, for an integer 
  $k \geq 2$. Assume furthermore that there are no divisors satisfying 
  the conditions $(**)$. Then 
\[ L \sim 2D + \Gamma, \]
  where 
  $\Gamma$ is a smooth curve satisfying $\Gamma ^2 = -2$ and $\Gamma.D =1$.
\end{lemma}

\begin{proof}
  We calculate $(L-2D)^2 =-2$. From Riemann-Roch and the fact that 
  $L-2D \geq 0$, we get that $L-2D \sim \Delta$ for some effective divisor 
  $\Delta$ such that $\Delta ^2 = -2$. We calculate $\Delta.D =1$ and 
  $\Delta.L =0$. Since $L$ is nef, $\Delta$ has to be supported on a union 
  of smooth rational curves.

  If there exists a smooth rational curve $\Gamma$, occurring as a component 
  of $\Delta$, such that $\Gamma.D =1$ (and necessarily $\Gamma.L =0$), we get 
\[ L^2 (2D + \Gamma)^2 = (4k+2)^2 =  (L.(2D + \Gamma))^2,\]
  whence $L \sim 2D + \Gamma$ by Hodge index theorem. 

  So it is sufficient to prove the existence of such a curve $\Gamma$. 

  Since $\Delta.D =1$, there has to exist a 
  smooth rational curve $\Gamma$ in $\Delta$ satisfying $\Gamma.D \geq 1$. If 
  $\Gamma.D \geq 2$, define $D':= D + \Gamma$. Then we have 
  $D'.L =D.L \leq  \frac{1}{2} L^2$, whence by Hodge index 
  $2{D'}^2 \leq L.D'$. We also get ${D'}^2 \geq  D^2 + 2$, 
  whence $D'.L < {D'}^2+ k_0 +1$. By the 
  assumptions that there are no divisors $D$ satisfying $(**)$, we must have 
  $D.\Gamma =1$, and we are done.
\end{proof}

\begin{prop} \label{negmain}
  Let $L$ be a globally generated, big line bundle on a $K3$ surface 
  $S$ satisfying $L^2 \geq 4k_0$ for an integer $k_0 \geq 0$. If there are no 
  divisors $D$ satisfying the conditions $(**)$ for $k=k_0$, then any 
  $0$-dimensional subscheme of $S$ 
  where the $k_0$-very ampleness of $L$ fails, has some point in common with 
  a smooth rational curve $\Gamma$ such that $\Gamma.L \leq k_0+1$. Hence 
  $L$ is birationally $k_0$-very ample. 
\end{prop}

\begin{proof}
  If  $\Z$ is any $0$-dimensional subscheme of $S$ of length $l+1$, 
  for an integer $l \leq k_0$, such that the map 
\[ H^0 (L) \hpil H^0 (L  \* \O _{\Z}) \]
  is not onto, but for any proper subscheme $\Z'$ of $\Z$, the map 
\[ H^0 (L) \hpil H^0 (L  \* \O _{\Z'}) \]
  is onto, then by the results in Section \ref{numcond} again, there exists an 
  effective divisor $D$ containing $\Z$ and satisfying the numerical 
  conditions $(*)$ for $k=l$, $D^2 \geq -2$ and $h^1 (D)=0$. By our 
  assumptions, we must either have $D^2 = -2$ or $(L^2,D^2) = (4k_0 +2, k_0)$
  and $l=k_0$. 
  In the first case, we are done as in the proof of Lemma \ref{neg}. So we 
  must treat the second case. We will actually show that in this case $\Z$ 
  will meet a smooth rational curve $\Gamma$ such that $\Gamma.L =0$. 

  By Lemma \ref{numlem}, iv), we have $L-2D \geq 0$, and by our 
  assumptions on $L^2$ and $D^2$, we must have $L.D = D^2 + k_0 + 1$, whence 
  $A= \emptyset$ in Lemma \ref{numlem} and the exact sequence (\ref{enes}). 

  By Lemma \ref{help1}, we have $L \sim 2D + \Gamma$, where 
  $\Gamma$ is a smooth rational curve such that $\Gamma.D =1$. Define 
  $F:= L-D \sim D + \Gamma$.

  Since $\Gamma.F = -1$, $\Gamma$ must be a base component of $F$. From this 
  it follows that $h^0 (F) = h^0 (D)$, and since $h^1 (D)=0$ and one 
  calculates $F^2 = D^2$, we get from Riemann-Roch that $h^1 (F) =0$.

  Tensorising the exact sequences (\ref{BSes}) 
  and (\ref{enes}) by $\O_X (-F)$ and $\O_X (-D)$, respectively, and using 
  $H^1 (\Gamma) = H^1 (F) =0$ and Serre duality, we find 
  $h^0 (\O_X (D) \* \I _ {\Z}) =1$ and $h^0 (\O_X (F) \* \I _ {\Z}) =2$ 
  respectively.

  The latter equality means that we can choose two distinct elements $F_1$ 
  and $F_2$ in $|F|$ both containing $\Z$ (scheme-theoretically). But since 
  $\Gamma$ is a base component of $|F|$, we must have 
  $F_1 =D_1 + \Gamma$ and  $F_2 =D_2 + \Gamma$, for two distinct elements 
  $D_1$  and $D_2$ of $|D|$. If $\Z$ meets $\Gamma$, we are done. If not, 
  we would have both $D_1$  and $D_2$ containing $\Z$ (scheme-theoretically). 
  But this contradicts the fact that $h^0 (\O_X (D) \* \I _ {\Z}) =1$.
\end{proof}

This concludes the proof of the equivalence of parts a)-d) in Theorem 
\ref{birkvathm}.

\section{The Clifford Index and Gonality of Smooth Curves in $|L|$}
\label{cliff}

We briefly recall the definition and some properties of gonality and Clifford 
index of curves.
Let $C$ be a smooth irreducible curve of genus $g \geq 2$.
We denote by $g^r_d$ a linear system of dimension $r$ and degree $d$ and say 
that $C$ is $k$-{\it gonal} (and that $k$ is its {\it gonality}) if $C$
posesses a $g^1_k$ but no $g^1_{k-1}$. In particular, we call a $2$-gonal 
curve {\it hyperelliptic} and a $3$-gonal curve {\it trigonal}. We denote by $\gon C$ the gonality of $C$. Note that if 
$C$ is $k$-gonal, all $g^1_k$'s must necessarily be base point free and 
complete.

If $A$ is a line bundle
on $C$, then the {\it Clifford index} of $A$ is the integer
\[ \Cliff A = \deg A - 2(h^0 (A) -1). \]
The {\it Clifford index of $C$} itself is defined as 
\[ \Cliff C = \min \{ \Cliff A \hspace{.05in} | \hspace{.05in} h^0 (A) \geq 2, h^1 (A) \geq 2 \}. \]
Clifford's theorem then states that $\Cliff C \geq 0$ with equality if
and only if $C$ is hyperelliptic and $\Cliff C =1$ if
and only if $C$ is trigonal or a smooth plane quintic.

At the other extreme, we get from Brill-Noether theory (cf.
\cite[V]{acgh}) that the gonality of $C$ satisfies $\gon C \leq \lfloor 
  \frac{g+3}{2} \rfloor $, whence $\Cliff C \leq \lfloor 
  \frac{g-1}{2} \rfloor $. For the general curve of genus
  $g$, we have $\Cliff C = \lfloor \frac{g-1}{2} \rfloor $.

We say that a line bundle $A$ on $C$ {\it contributes to the Clifford
  index of $C$} if $h^0 (A), h^1 (A) \geq 2$ and that it {\it computes 
the Clifford index of $C$} if in addition $\Cliff C = \Cliff A$. 

Note that $\Cliff A = \Cliff \omega _C \* A^{-1}$.

The {\it Clifford dimension} of $C$ is defined as
\[ \min \{ h^0 (A) -1 \mbox{  } | \mbox{  } A  
\mbox{ computes the Clifford index of } C \} . \]

A line bundle $A$ which achieves the minimum and computes the Clifford index, is said to {\it compute} the Clifford dimension. A curve of Clifford index $c$ is $(c+2)$-gonal if and only if it has Clifford dimension $1$. For a general curve $C$, we have  $\gon C = c+2$.

\begin{lemma} \label{gon1}
  \cite[Theorem 2.3]{cm} The gonality $k$ of a smooth irreducible projective curve $C$ satisfies 
\[ \Cliff C +2 \leq k \leq \Cliff C +3.\]
\end{lemma}

The curves satisfying $\gon C = \Cliff C +3$ are conjectured to be very rare and called {\it exceptional} (cf. \cite[(4.1)]{ma}).

\begin{lemma} \label{gon2}
  \cite[Prop. 1.2]{mp} A smooth irreducible projective curve $C$ is $k$-gonal if and only if $K_C$ is $(k-2)$-very ample but not $(k-1)$-very ample.
\end{lemma}

Recall also the result of Green and Lazarsfeld \cite{gl}, which states
that if $L$ is a base point free line bundle on a $K3$ surface $S$, then
$\Cliff C$ is constant for all smooth irreducible $C \in |L|$, and if
$\Cliff C < \lfloor \frac{g-1}{2} \rfloor $, then there exists a line
bundle $M$ on $S$ such that $M_C := M \* \O _C$ computes the Clifford index of
$C$ for all smooth irreducible $C \in |L|$. (Note that since $(L-M) \*
\O _C \iso \omega _C \* {M_C}^{-1}$, the result is symmetric in $M$ and
$L-M$.)

It turns out that we can choose $M$ so that it satisfies certain
properties. We will need the following result in the sequel.

\begin{lemma} \label{chooseD}
  Let $L$ be a  base point free line bundle on a $K3$ surface $S$ with
  $L^2 = 2g-2 \geq 2$. Let $c$ be the Clifford index of any smooth
  curve in $|L|$.  

  If $c < \lfloor \frac{g-1}{2} \rfloor $, then there exists a smooth
  curve $D$ on $S$ satisfying $0 \leq D^2 \leq c+2$, $2D^2 \leq D.L$
  (either of the latter two inequalities being an equality if and only if 
  $L \sim 2D$) and
\[ \Cliff C = \Cliff (\O _S (D) \* \O_C) = D.L - D^2 -2 \] 
  for any smooth curve $C \in |L|$. 
\end{lemma}

\begin{proof}
  It follows from the proof of the main theorem in \cite{gl}, as worked
  out by Martens in \cite[(2.3)]{ma}, that we can choose the line bundle
  $M$ above so that $M \* \O_C$ and $(L-M) \* \O_C$ compute 
  the Clifford index of $C$ and furthermore $h^0 (M) = h^0 (M_C) \geq 2$, $h^0
  (L-M)= h^0 (\omega_C \* M_C ^{-1}) =  h^1 (M_C) \geq 2$, 
  $h^1 (M) = h^1(L-M)=0$ and such    
  that either $|M|$ or $|L-M|$ is base point free (any one of them by our
  choice, but not necessarily both at the same time, unless $L$ is ample
  \footnote{Here in fact there is a slight malformulation in Cor. 2.3 in
  \cite{ma}, where it seems as if one can choose both $|M|$ and $|L-M|$ to
  be base point free at the same time, which is not the case.}).

  By symmetry, we can assume $M.L \leq (L-M).L$, or equivalently $2M.L
  \leq L^2$,  and we will choose $M$ 
  to be base point free. Since $h^1 (M)=0$, it follows from Proposition \ref{sd1}
  again that the generic member of $|M|$ is a smooth irreducible curve
  $D$ 
  of genus $\geq 0$. Combining $2D.L \leq L^2$ with the Hodge index
  theorem, we get $2D^2 \leq D.L$, with equality if and only if $L \sim 2D$.

  Furthermore, by the exact cohomology sequence associated to the
  standard short exact sequence 
\[ 0 \hpil \O _S (D-C) \hpil \O _S (D) \hpil \O _S (D) \* \O_C  \hpil 0,
  \]
  and using that $h^1 (D-C) = h^1 (L-M) = 0$, we get from Riemann-Roch
\[ h^0 (\O _S (D) \* \O_C) = h^0 (D) = \frac{1}{2}D^2 +2, \]
  so 
\[ \Cliff C = \deg (\O _S (D) \* \O_C) - 2( h^0 (D)-1) = D.L - D^2 -2. \]

Since $2D^2 \leq D.L$ we finally get $c = D.L - D^2 -2 \geq D^2-2$. 
\end{proof}

\begin{rem} \label{remcliff}
  {\rm From the exact sequence above it follows that if $D$ is any divisor
  such that $h^0 (D),h^0 (L-D) >0$, then $h^0 (\O _S (D) \* \O_C) \geq h^0
  (D)$, $h^1 (\O _S (D) \* \O_C) \geq h^0 (L-D)$ and } 
  $\Cliff (\O _S (D) \* \O_C) \leq D.L - D^2 -2$. 
\end{rem}

We have the following result

\begin{lemma} \label{cliffeq}
  Let $L$ be a globally generated and big line bundle on a $K3$ surface 
  and denote by $c$ the Clifford index of any smooth curve in $|L|$. Let 
  $k_1$ be the smallest integer such that the conditions $(*)$ are fulfilled 
  for an integer $k \geq 1$ and a divisor $D$ satisfying $D^2 \geq 0$, and 
  $k_2$ be the smallest integer such that $L^2 < 4k_2$. Define 
  $k_0 := \min \{ k_1,k_2 \}$

 Then $c=k_0-1$.
\end{lemma}

\begin{proof}
  Let $g=\frac{1}{2}L^2+1$ be the genus
  of all the smooth curves in $|L|$.

  If $L^2 < 4k_2$, then we have 
\[ \lfloor \frac{g-1}{2} \rfloor = \lfloor \frac{L^2}{4} \rfloor < k_2,\]
  whence $c \leq k_2-1$.
  
  If $k_1$ is the smallest integer such that the conditions $(*)$ are 
  fulfilled for an integer $k \geq 1$ and a divisor $D$ satisfying 
  $D^2 \geq 0$, then by Proposition \ref{smcurprop} we can find a smooth
  curve $D_0$ satisfying $(*)$ for $k=k_1$ with equality in the
  middle and 
  such that $h^0(F) \geq h^0(D_0) \geq 2$ and $h^1(D_0)=0$, where 
  $F := L-D_0$. Let $C$ be any smooth curve in $|L|$.
  From Remark \ref{remcliff} we have
  \[ \Cliff (\O_C (D_0)) \leq  D_0.C - D_0^2 -2 = k_1-1, \]
  whence $c \leq k_1-1$ again.

  We now prove the opposite inequality.

  If $c= \lfloor \frac{g-1}{2} \rfloor$, we have $L^2 = 4c$ or $4c+2$,
  whence $k_2 = c+1$.

  If $c < \lfloor \frac{g-1}{2} \rfloor$, there exists by Lemma
  \ref{chooseD} a smooth
  curve $D$ satisfying $0 \leq D^2 \leq c+2$, $2D^2 \leq D.L$ with
  equality if and only if $L \sim 2D$, and also 
  $c= D.L - D^2 -2$. But this means that $D$ satisfies the
  conditions $(*)$ for $k=c+1$, whence $k_1 \leq c+1$.
\end{proof}

The following proposition investigates the particular case 
$(L^2, D^2) = (4k+2,k)$ apperaring in the conditions $(**)$.

\begin{prop} \label{onlyex}
  Let $L$ be a spanned and big line bundle on a $K3$ surface and denote by 
  $c$ the Clifford index of all smooth curves in $|L|$. Then the following 
  two conditions are equivalent:
           \begin{itemize}
  \item [i)]  $L \sim 2D + \Gamma$, 
              where $D$ is an effective divisor and 
              $\Gamma$ is a smooth curve satisfying $D^2 = k$, 
              $\Gamma ^2 = -2$ and $\Gamma.D =1$, for an integer $k \geq 1$
              and there are no divisors satisfying the conditions $(**)$,
  \item [ii)] $c = k-1$ and all smooth curves in $|L|$ have gonality $c+3$ 
              (equivalently, all smooth curves in $|L|$ have Clifford 
              dimension $>1$).
           \end{itemize}
\end{prop}

\begin{proof}
  Given i), by Remark \ref{remcliff} we have $c \leq D.L - D^2 -2 =k-1$. Since 
  there are no divisors satisfying the conditions $(**)$, we must have equality.
  
  By Lemma \ref{gon1}, we have $k+1 \leq \gon C \leq k+2$. We saw in the 
  proof of Proposition \ref{negmain} that if $\Z$ is any $0$-dimensional scheme where the $k$-spannedness of $L$ fails, then $\Z$ has a point in common with $\Gamma$. But $\Gamma.L=0$, whence $\Gamma$ does not intersect any smooth curve in $|L|$ and $L_{|C} \iso \omega_C$ is $(k+1)$-spanned for all smooth $C \in |L|$. Therefore $\gon C \geq k+2$ by Lemma \ref{gon2}.

 Conversely, given ii), the curves in $|L|$ cannot have general Clifford 
 index, so $c=k-1 < \lfloor \frac{L^2}{4} \rfloor $, whence $L^2 \geq 4k$ and by Lemma \ref{cliffeq} there exists a divisor and, a posteriori, a smooth curve $D$ satisfying the conditions $(*)$ and $D^2 \geq 0$. 

Assume $D$ satisfies the conditions $(**)$ as well (i.e. 
$(L^2, D^2) \not = (4k+2,k)$). In the proof of Lemma \ref{smcurlem} we proved 
that the line bundle $F_D$ is base point free on $D$, where $F:=L-D$ and 
$F_D:=F \* \O_D$ as usual. Furthermore, since $F$ computes the Clifford index 
of any smooth curve in 
$|L|$, one sees that $h^1 (F) =0$ and that the possibly empty base locus 
$\Delta$ of $|F|$ must satisfy $\Delta.L =0$. Hence, by \cite[Lemma 2.2]{cp}, 
we have that there exists a smooth curve in $|L|$ of gonality 
$F.D = k+1 = c+2$, a contradiction.

Therefore, $(L^2, D^2) = (4k+2,k)$ and there are no divisors satisfying 
$(**)$. 
By Lemma \ref{help1}, $L$ has the desired form.
\end{proof}

\begin{rem} \label{eisex}
  {\rm The example above is essentially the same as the one treated in 
\cite[Thm. 4.3]{elms}, where the authors prove that all smooth curves in 
$|L|$ are exceptional under the stronger hypothesis that 
$\Pic S \iso \ZZ D \+ \ZZ \Gamma$. We have proved this under the weaker 
hypothesis that there are no divisors satisfying 
$(**)$ and also proved  that this is the only example of a complete linear 
system on a $K3$ surface where all smooth curves are exceptional. 

It should be noted, however, that the authors in \cite{elms} also prove that 
the Clifford dimension of all smooth curves in $|L|$ is 
$\frac{1}{2} D^2 +1= \frac{1}{2}k +1$, a statement which is absent in our 
treatment. } 
\end{rem}


Now the equivalence of parts e)-g) with a)-d) in Theorem \ref{birkvathm} follows rather easily.

It is clear from Lemma \ref{gon2} that the conditions e) and g) are 
equivalent. The equivalence between d) and f) follows from Lemma \ref{cliffeq} and Proposition \ref{onlyex}, and the equivalence between e) and f) follows from Proposition \ref{onlyex}.

What remains now is condition h). There are two approaches here: one is to study the vector bundles $E$ appearing in Section \ref{numcond}, the other is to study vector bundles naturally arising from the study of line bundles on smooth curves on $K3$ surfaces. We will follow the second approach. For details we refer to \cite{la}, \cite{gl} and \cite{cp}.

Let $C$ be a smooth curve on a $K3$ surface $S$. Recall that if $A$ is a line bundle on $C$ with the property that both $A$ and $\omega_C \* A^{-1}$ are generated by their global sections, then one can associate to the pair $(C,A)$ a vector bundle $E(C,A)$ of rank $h^0 (A)$ as follows. Thinking of $A$ as a coherent sheaf on $S$, we get a short exact sequence 
\[ 0 \khpil F(C,A) \khpil H^0 (A) \* _{\CC} \O_S  \khpil A \khpil 0 \]
of $\O_S$-modules, where $F(C,A)$ is locally free (since $A$ is locally isomorphic to $\O_C$ and hence has homological dimension one over $\O_S$).

The vector bundle $E(C,A) := {F(C,A)}^*$ will have rank $h^0 (A)$ and the following properties: $H^1(E(C,A)) =  H^2(E(C,A))=0$, $\det E(C,A) = L$, $c_2 (E(C,A))= \deg A$ and $E(C,A)$ is generated by its global sections.

Note that if $A$ is any line bundle on $C$ computing the Clifford index or the gonality of $C$, then $A$ will satisfy the conditions that  both $A$ and $\omega_C \* A^{-1}$ are generated by their global sections, and we can carry out the construction of the vector bundle $E(C,A)$ above. 

Also recall that if $E$ is any vector bundle on $S$ generated by its global sections and satisfying  $H^1(E) =  H^2(E)=0$, then $E= E(C,A)$ for some pairs $(C,A)$, where $C$ is a smooth curve in $|\det E|$ and $A$ is a line bundle on $C$ satisfying $h^0 (A) = \rank E$ and $\deg A = c_2 (E)$ \cite[Lemma 1.2]{cp}.

Now the equivalence between e) and h) in Theorem \ref{birkvathm} follows immediately. Indeed, the existence of a smooth curve in $|L|$ of gonality $ \leq k+1$ is equivalent to the existence of a rank two vector bundle $E$ generated by its global sections and satisfying $H^1(E) =  H^2(E)=0$, $\det E = L$ and $c_2 (E)  \leq k+1$.

This concludes the proof of Theorem \ref{birkvathm}.

\noindent \textit{Mathematics Subject Classification:} 14J28.

\noindent \textit{Key words:} Line bundles, $K3$ surfaces, Enriques surfaces,
$k$-very ampleness, $k$-spannedness, Clifford index, gonality. 

\noindent \textit{E-mail:} andreask@mi.uib.no.

\end{document}